\documentclass{amsart}

\usepackage{amssymb}
\usepackage{amsmath}
\usepackage{amsfonts}
\usepackage{latexsym}
\usepackage{amsthm}
\usepackage{enumerate}

 \newtheorem{thm}{Theorem}[subsection]
 \newtheorem{thm*}{Theorem}
 \newtheorem{cor}[thm]{Corollary}
 \newtheorem{lem}[thm]{Lemma}
 \newtheorem{prop}[thm]{Proposition}
 \theoremstyle{definition}
 \newtheorem{defn}[thm]{Definition}
 \theoremstyle{remark}
 
 \newtheorem{exam}[thm]{Example}
 \numberwithin{equation}{subsection}

 \DeclareMathOperator{\RE}{Re}
 \DeclareMathOperator{\IM}{Im}

 \newcommand{\R}{\ensuremath{\mathbb{R}}}
 \newcommand{\C}{\ensuremath{\mathbb{C}}}
 \newcommand{\D}{\ensuremath{\mathbb{D}}}
 \newcommand{\Q}{\ensuremath{\mathbb{Q}}}
 \newcommand{\Z}{\ensuremath{\mathbb{Z}}}

 \newcommand{\pp}{\ensuremath{\mathbb{P}}}

 \newcommand{\derp}[2]{\ensuremath{\frac{\partial #1}{\partial #2}}}
 \newcommand{\derpp}[1]{\ensuremath{\derp{}{#1}}}
\begin{document}

\title[nacs and non-K\"{a}hler compact complex manifolds]
 {Normal almost contact structures and non-K\"{a}hler compact complex manifolds}

\author{ M\`{o}nica Manjar\'{\i}n}

\address{Department d'\`{A}lgebra i Geometria, Facultat de Matem\`{a}tiques, Universitat de Barcelona,
Gran Via de les Corts Catalanes 585 Barcelona 08007, Spain}

\email{manjarin@ub.edu}




\date{December 18, 2006.}

\thanks{}


\begin{abstract}
We construct some families of complex structures on compact
manifolds by means of normal almost contact structures (nacs) so
that each complex manifold in the family has a non-singular
holomorphic flow.  These families include as particular cases the
Hopf and Calabi-Eckmann manifolds and the complex structures on the
product of two normal almost contact manifolds constructed by
Morimoto. We prove that every compact K\"{a}hler manifold admitting
a non-vanishing holomorphic vector field belongs to one of these
families and is a complexificacion of a normal almost contact
manifold. Finally we show that if a complex manifold obtained by our
constructions is K\"{a}hlerian the Euler class of the nacs (a
cohomological invariant associated to the structure) is zero. Under
extra hypothesis we give necessary and sufficient conditions for the
complex manifolds so obtained to be K\"{a}hlerian.
\end{abstract}

\maketitle

\section{Introduction}

Most of the known examples of complex manifolds, in particular all
projective manifolds, are of K\"{a}hler type. Nevertheless, the
existence of a K\"{a}hler metric imposes strong topological
restrictions on the manifold, for instance its odd Betti numbers are
even. Riemann surfaces are always K\"{a}hlerian and compact complex
surfaces if and only its first Betti number is even (c.f.
\cite{Buch},\cite{Lama}). For higher dimensions there is not a
simple characterization of K\"{a}hler manifolds, however one would
expect that they are rather the exception than the rule. For
instance a corollary of a result by Taubes (c.f \cite{Tau}) implies
that every finite presentation group is the fundamental group of a
non-K\"{a}hler compact complex 3-manifold. Historically the first
examples of non-K\"{a}hler manifolds were constructed by H. Hopf as
a quocient of $\C^n-\{0\}$ for $n>1$ by a contracting holomorphism
of $\C^n$ which fixes the origin. Later, E. Calabi and B. Eckmann
described a class of non-K\"{a}hler complex structures on the
product $S^{2n+1}\times S^{2m+1}$ for $n,m\geq 0$ such that the
corresponding complex manifold is the total space of a holomorphic
elliptic principal bundle over $\pp^n\times \pp^m$. In
\cite{MNLoeb1} J.J. Loeb and M. Nicolau generalized Calabi-Eckmann
and Hopf structures by the construction of a class of complex
structures on the product $S^{2n+1}\times S^{2m+1}$ that contains
the precedents. Similar techniques have been used by S.L\'{o}pez de
Medrano and A.Verjovsky in \cite{LMVer} to construct another family
of non-K\"{a}hlerian compact manifolds and later generalized by
L.Meersseman in \cite{Meer}.

The second section of the paper is devoted to describe a general
procedure to obtain complex manifolds by means of elementary
geometrical constructions. We depart from odd-dimensional manifolds
admitting a normal almost contact structure (nacs for shortness),
i.e. a CR-structure of maximal dimension and a transverse CR-action
of $\R$. More precisely we consider three cases: (\textbf{A})
products of two real manifolds endowed with a nacs, (\textbf{B})
$S^1$-principal bundles over a manifold with a nacs (with an extra
restriction on the bundle) and (\textbf{C}) suspensions of a
manifold with a nacs by a suitable automorphism. In particular we
generalize Morimoto's construction of a complex structure on a
product of two normal almost contact manifolds (c.f. \cite{Morimo}).
The constructions of cases \textbf{A} and \textbf{C} produce a
compact complexification of the original normal almost contact
manifold $\mathrm{M}$, i.e. a compact complex manifold $\mathrm{X}$
such that $\mathrm{M}$ is a real submanifold of $\mathrm{X}$ so that
its CR-structure is compatible with the complex structure of
$\mathrm{X}$ and a holomorphic vector field on $\mathrm{X}$ whose
real part is the vector field of the CR-action on $\mathrm{M}$.
Moreover we prove that given a compact K\"{a}hler mani\-fold
admitting a holomorphic vector field without zeros its complex
structure can be recovered by the construction of case \textbf{C}.
Therefore, every compact K\"{a}hler manifold admitting a
non-vanishing holomorphic vector field is a compact complexification
of a nacs. We also show that double suspensions of compact complex
manifolds by two commuting automorphisms (see p.8), which can be
obtained by means of the construction of case \textbf{C}, present a
remarkable property. Namely, every compact K\"{a}hler manifold
admitting a non-vanishing holomorphic vector field can be endowed
with a complex structure on the underlying smooth manifold
arbitrarily close to the original one which turns it into a double
suspension.

In the third section we study criteria to determine when the complex
manifolds of the above three families are K\"{a}hlerian in terms of
properties of the departing nacs. The common feature of all the
complex manifolds that we construct is the existence of a
holomorphic vector field without zeros. Using this fact we will show
that for these complex structures to be K\"{a}hlerian the Euler
class of the nacs must be zero. When the flows associated to the
nacs are isometric we prove that a compact complex manifold obtained
by the constructions of cases \textbf{A} or \textbf{B} is K\"{a}hler
if and only if the Euler class (or classes) is zero and the flows
(or flows) is transversely K\"{a}hlerian. For suspensions (case
\textbf{C}) of a normal almost contact manifold $\mathrm{M}$ by an
automorphism $f$ we give a complete characterization when the
CR-structure is Levi-flat and $f^{\ast}=id$ acting on
$H^1(\mathrm{M},\C)$. Finally, we prove that a double suspension is
K\"{a}hlerian if and only if the departing manifold is K\"{a}hlerian
and the two automorphisms preserve a K\"{a}hler class.

Throughout the paper all manifolds are supposed to be smooth and
connected and all differentiable objects to be of class
$\mathcal{C}^{\infty}$.

I would like to thank my advisor, Marcel Nicolau, for having
proposed me this problem, for having guided me during these years
and for many useful comments regarding this paper. I am also
grateful to Marco Brunella for suggesting the characterization for
double suspensions and to Aziz El Kacimi for pointing out its
validity in a more general context.

\section{Compact complex structures defined from nacs}

\subsection{Normal almost contact structures
(nacs)}\label{section:nacs}

Recall that a complex subbundle $\Phi^{1,0}$ of dimension $m$ of the
complexified tangent bundle $T^{\C}\mathrm{M}=T\mathrm{M}\otimes \C$
of a manifold $\mathrm{M}$ is called a CR-\emph{structure} on M of
dimension $m$ (cf. \cite{KobNo2} or \cite{Bogg}) if: \pagebreak
\begin{enumerate}[\bf (i)]
\item $\Phi^{1,0}\cap \overline{\Phi^{1,0}}=\{0\}$;
\item $\Phi^{1,0}$ is involutive, i.e. $[\Phi^{1,0},\Phi^{1,0}]\subset
\Phi^{1,0}$.
\end{enumerate}
The complex bundle $\Phi^{1,0}$ induces a real subbundle
$\mathcal{D}=T\mathrm{M}\cap (\Phi^{1,0}\oplus
\overline{\Phi^{1,0}})$ of $T\mathrm{M}$. We define an endomorphism
$J:{\mathcal D}\rightarrow {\mathcal D}$ imposing that
$v-\mathrm{i}Jv\in \Phi^{1,0}$ for every $v\in {\mathcal D}$. Note
that we can determine the CR-structure by giving
$(\mathrm{M},{\mathcal D}, J)$. Setting
$\Phi^{0,1}=\overline{\Phi^{1,0}}$ we have a decomposition
${\mathcal D}\otimes \C=\Phi^{1,0}\oplus \Phi^{0,1}$ where
$\Phi^{1,0}$ and $\Phi^{0,1}$ are the eigenspaces of $J$ (extended
by complex linearity to ${\mathcal D}\otimes \C$) of eigenvalue
$\mathrm{i}$ and $-\mathrm{i}$ respectively. We denote by
$\mathrm{Aut}_{\mathrm{CR}}(\mathrm{M})$ the subset of
$\mathrm{Diff(M)}$ of maps $f$ such that $df$ preserves ${\mathcal
D}$ and commutes with $J$. Let $\{\varphi_t: \,t\in \R\}$ be the
flow induced by a smooth $\R$-action on M. We say that
$\{\varphi_t\}$ defines a CR-\emph{action}\index{CR-action} if
$\varphi_t\in\mathrm{Aut}_{\mathrm{CR}}(\mathrm{M})$ for each $t$.
When $\dim_{\R}\mathrm{M}=2m+1$ we call the action \emph{transverse}
to the CR-structure if the smooth vector field
$T=(d\varphi_t/dt)_{t=0}$ is everywhere transverse to ${\mathcal
D}$, i.e. $\langle T \rangle \oplus {\mathcal D}$ has real dimension
$2m+1$ at every point.

\begin{defn}\label{defn:nacs}
A \emph{normal almost contact structure} (or \emph{nacs}) on a manifold $\mathrm{M}$
of odd-dimension is a pair $(\Phi^{1,0},\varphi_t)$ where $\Phi^{1,0}$ is a
CR-structure of maximal dimension and $\{\varphi_t\}$ a flow induced by a smooth
$\R$-action defining a transverse CR-action. Given a nacs $(\Phi^{1,0},\varphi_t)$
we define its \emph{characteristic 1-form} $\omega$ by the conditions $\omega(T)=1$
and $\ker \omega=\mathcal{D}$ (where $T=(d\varphi_t/dt)_{t=0}$ and
$\mathcal{D}=T\mathrm{M}\cap (\Phi^{1,0}\oplus \overline{\Phi^{1,0}})$) and its
\emph{associated flow} $\mathcal{F}$ as the flow induced by $\{\varphi_t\}$, which
is transversely holomorphic.
\end{defn}

Alternatively a nacs can be determined by an endomorphism $\varphi$ on the tangent
space, a vector field $T$ and a 1-form $\omega$. The tercet $(\varphi,T,\omega)$ is
called an almost contact structure on $\mathrm{M}$ if: \textbf{(1)} $\omega(T)=1$,
\textbf{(2)} rank $\varphi=2n$, \textbf{(3)} $\varphi(T)=0$, \textbf{(4)}
$\omega(\varphi(X))=0$ and $\varphi^2(X)=-X+\omega(X)T$ for every tangent vector
field $X$ on $\mathrm{M}$. There is an almost contact structure on $\R$ given
by~$(0,\derpp{t},dt)$. If $\mathrm{M}_1$ has an almost contact structure
$(\varphi_1,T_1,\omega_1)$ then there is an almost complex structure~$K$ on
$\mathrm{M}_1\times \R$ defined by
$K(X_1,a\derpp{t})=(\varphi_1(X_1)-aT_1,\omega_1(X_1)\derpp{t})$. The almost contact
structure on $\mathrm{M}_1$ is called \emph{normal} if $K$ is integrable (cf.
\cite{Blair}). It is not difficult to see that the two definitions are equivalent.

Recall that a form $\omega\in \Omega^{\ast}(\mathrm{M})$ is called
\emph{basic} with respect to a foliation $\mathcal{F}$ if
$i_S\omega=i_Sd\omega=0$ for every vector field $S$ tangent to the
leaves of $\mathcal{F}$.

\begin{lem}
Let ${\mathcal F}$ be a transversely holomorphic flow on a compact manifold
$\mathrm{M}$ generated by a real vector field $T$ without zeros and let $\omega$ be
a $1$-form such that $\omega(T)=1$. Set ${\mathcal D}=\ker \omega$ and let $J$ be
the almost-complex structure on ${\mathcal D}$ induced by ${\mathcal F}$. Then
$({\mathcal D},J,T)$ is a nacs if and only if $i_Td\omega=0$ and the basic form
$d\omega$ is of type $(1,1)$ with respect to the holomorphic structure transverse to
${\mathcal F}$.
\end{lem}

\begin{proof}
Let $\Phi^{1,0}$ the vectors in ${\mathcal D}^{\C}$ of type $(1,0)$
with respect to $J$. The vector field $T$ preserves ${\mathcal D}$,
i.e. $[T,{\mathcal D}]\subset {\mathcal D}$ if and only if $i_T
d\omega=0$. In this case $\Phi^{1,0}$ defines a CR-structure, i.e.
$[\Phi^{1,0},\Phi^{1,0}]\subset \Phi^{1,0}$, if and only if
$d\omega$ is of type $(1,1)$ with respect to the complex structure
transverse to ${\mathcal F}$. It is clear that then $T$ defines a
transverse CR-action.
\end{proof}

Let $\mathrm{M}$ be a compact manifold and ${\mathcal F}$ a transversely holomorphic
isometric flow on $\mathrm{M}$ defined by a Killing vector field $T$. If there
exists a characteristic 1-form $\omega$ (which verifies $\omega(T)=1$ and
$i_Td\omega=0$) such that $d\omega$ is of type (1,1) then $\mathrm{M}$ admits a
nacs. Analogously a $S^1$-principal bundle over a compact complex manifold admits a
nacs provided that we can choose a connection 1-form such that its curvature form is
of type $(1,1)$. It is well known that a $S^1$-principal bundle admits such a
connection form if and only if it is the unit bundle associated to a hermitian
metric on a holomorphic line bundle. More generally, compact Seifert fibrations over
a complex orbifold also provide examples of transversely holomorphic isometric
flows, therefore they admit a nacs provided that there exists a suitable
characteristic 1-form. Notice that if $\dim_{\R}\mathrm{M}=3$ the last condition of
the lemma is always fulfilled since every 2-form on a compact Riemann surface is of
type $(1,1)$.

It is also known that a compact connected Lie group $\mathrm{K}$ of
odd dimension greater than one always admits a non left-invariant
nacs (cf. \cite{LMN}).

\medskip

Let now $({\mathcal D},J)$ be a CR-structure on $\mathrm{M}$ and
suppose that the distribution ${\mathcal D}$ is a contact structure,
i.e. the characteristic 1-form $\omega$ verifies $\omega\wedge
(d\omega)^n\neq 0$. Then $({\mathcal D},J)$ is called a
\emph{strictly pseudo-convex} CR-\emph{structure} on M. The couple
of a strictly pseudo-convex CR-structure of maximal dimension and a
transverse CR-action of $\R$ on an odd-dimensional manifold is also
known as a \emph{normal contact structure}. For compact connected
3-manifolds it is known that if they admit a normal contact
structure the vector field defining the CR-action is Killing (c.f.
\cite{Blair}). The opposite situation to a strictly pseudo-convex
CR-structure from the point of view of the real integrability of the
distribution ${\mathcal D}$ is Levi-flatness, that is, the condition
$\beta \wedge d\beta \equiv 0$ where $\beta$ is a 1-form such that
$\ker\beta=\mathcal{D}$ or equivalently $d\omega=0$ where $\omega$
is the characteristic 1-form of the nacs. In this case we can easily
construct examples of nacs such that its associated flow is not
isometric. Recall that the \emph{suspension} of a compact manifold
$\mathrm{N}$ by $g\in\mathrm{Diff}(\mathrm{N})$ is the compact
manifold $\mathrm{N}\times_g\R$ given by $\mathrm{N}\times \R/\sim$
where $(z,s)\sim (g(z),s+1)$. When $\mathrm{N}$ is a compact complex
manifold and $g\in Aut_{\C}(\mathrm{N})$ the suspension
$\mathrm{N}\times_g\R$ carries a natural Levi-flat CR-structure
defined by $T\mathrm{N}$ and a transverse CR-action induced by
$\derpp{s}$. If we choose $g$ such that it is not an isometry for
any metric on N, for instance $\mathrm{N}=\C P^1$ and $g(z)=\lambda
z$ with $\lambda\in\C$ such that $|\lambda|\neq 1$, the flow
${\mathcal F}$ generated by $\derpp{s}$ is clearly not isometric for
any Riemannian metric on $\mathrm{N}\times_g\R$.

\subsection{The Euler Class}

We introduce here the notion of Euler class, which generalizes the
classical notion of Euler class of an isometric flow.

\begin{defn}\label{def:EulerClass}
Let M be a compact manifold endowed with non-vanishing vector field
$T$ and a 1-form $\omega$ such that $\omega(T)=1$ and
$i_Td\omega=0$. We denote by $\mathcal{F}$ the flow induced by $T$.
We define the \emph{Euler class} of the pair $(\mathrm{M},T)$ as the
basic cohomology class given by
$$e_{\mathcal F}(\mathrm{M})=[d\omega]\in
H^2(\mathrm{M}/\mathcal{F}).$$
\end{defn}

Recall that the \emph{basic cohomology}
$H^{\ast}(\mathrm{M}/\mathcal{F})$ is defined as the cohomology of
the differential complex $(\Omega^{\ast}(\mathrm{M}/\mathcal{F}),d)$
of \emph{basic forms} for $\mathcal{F}$. Note that the class
$e_{\mathcal F}(\mathrm{M})$ depends on the vector field $T$ but not
on the 1-form, provided that it verifies the above conditions. In
particular we can consider the Euler class of a nacs. As the
vanishing of the Euler class will be a necessary condition for the
complex manifolds that we construct to be K\"{a}hlerian we will
discuss some criteria to determine when it is zero.

\begin{lem}
Let $\mathrm{M}$ be a compact manifold endowed with non-vanishing
vector field $T$ and a 1-form $\omega$ such that $\omega(T)=1$ and
$i_Td\omega=0$. We denote by $\mathcal{F}$ the flow induced by $T$.
Then the following conditions are equivalent:\pagebreak
\begin{enumerate}[\bf (a)]
\item $e_{\mathcal F}(\mathrm{M})=0$.
\item There exists a \emph{closed} 1-form $\chi$ on $\mathrm{M}$ such that
$\chi(T)=1$.
\item There exists a distribution transverse to $T$ of maximal dimension and invariant by
the flow which is \emph{integrable}.
\end{enumerate}
\end{lem}

\begin{proof}
To prove $\mathbf{(a)\Rightarrow (b)}$ note that if
$e_{\mathcal{F}}(\mathrm{M})=0$ then there exists a basic 1-form
$\alpha$ tal que $d\alpha=d\omega$. Then $\chi=\omega-\alpha$ is a
closed 1-form such that $\chi(T)=1$. To see that
$\mathbf{(b)\Rightarrow (c)}$ it is enough to remark that the
distribution is given by $\ker \chi$ and since $L_T \chi=0$ the form
$\chi$ is invariant by the flow. To prove $\mathbf{(c)\Rightarrow
(a)}$ one defines the 1-form $\chi$ imposing that it vanishes on the
distribution and $\chi(T)=1$. As $d\chi=0$ we have
$e_{\mathcal{F}}(\mathrm{M})=0$.
\end{proof}

\begin{cor}
Let $\mathrm{M}$ be a compact manifold endowed with a nacs with Levi-flat
\emph{CR}-structure. Then $e_{\mathcal{F}}(\mathrm{M})=0$.
\end{cor}

\begin{prop}\label{cor:tischler} Let $\mathrm{M}$ be a compact manifold endowed with a nacs.
If $\, e_{\mathcal F}(\mathrm{M})=0$ then $\mathrm{M}$ is a fiber
bun\-dle over $S^1$. In particular $b_1(\mathrm{M})\neq 0$ and
$\mathrm{M}$ is not simply connected.
\end{prop}

\begin{proof}
The first statement is a consequence of a theorem by D. Tischler
\cite{Tisch} that states that if $\mathrm{M}$ is a compact manifold
admitting a non-vanishing closed 1-form then $\mathrm{M}$ is a fibre
bundle over $S^1$. The second statement is an immediate consequence
of the homotopy exact sequence associated to a fibration.
\end{proof}

\begin{prop}\label{prop:contacte} Let $\mathrm{M}$ be a compact manifold endowed with a nacs
and $T$ the vector field inducing the \emph{CR}-action. If there exists a contact
form $\chi$ on $\mathrm{M}$ such that $\chi(T)=1$ and $i_Td\chi=0$ then $e_{\mathcal
F}(\mathrm{M})\neq 0$.
\end{prop}

The proof is analogous to the one of the corresponding statement for
isometry flows due to Saralegui (cf. \cite{Sara}).

\begin{cor}\label{cor:contacte}
Let $\mathrm{M}$ be a compact manifold endowed with a normal contact
structure. Then $e_{\mathcal{F}}(\mathrm{M})\neq 0$.
\end{cor}

\subsection{Geometrical constructions of complex structures}

\begin{prop}\label{prop:caseA}\emph{\textbf{(Case A)}} Let $\mathrm{M}_1$ and $\mathrm{M}_2$
be two manifolds endowed with a nacs. There exists a 1-parametric
family of complex structures $K_{\tau}$ on the product
$\mathrm{M}_1\times \mathrm{M}_2$ for $\tau\in\C\backslash\R$ so
that the complex manifold $\mathrm{M}_1\times \mathrm{M}_2$ admits a
non-vanishing holomorphic vector field $v$.
\end{prop}

\begin{proof}
Let us denote by $(T_1,\omega_1)$ and $(T_2,\omega_2)$ the vector
fields and the characteristic 1-forms of the nacs on $\mathrm{M}_1$
and $\mathrm{M}_2$ respectively. The 2-foliation ${\mathcal F}$ on
$\mathrm{M}_1\times\mathrm{M}_2$ generated by $T_1$ and $T_2$ is
transversely holomorphic and $[T_1,T_2]=0$. We set a distribution
${\mathcal D}$ given by $\ker\omega_1 \oplus \ker\omega_2$ and we
define a complex-valued 1-form
$\chi=\frac{i}{2\IM\tau}(\bar{\tau}\omega_1+\omega_2)$. The 2-form
$d\chi$ is basic with respect to $\mathcal{F}$ and $d\chi$ is of
type $(1,1)$ with respect to the transverse holomorphic structure of
$\mathcal{F}$. We define an almost complex structure
$\mathrm{K}_{\tau}$ on $\mathrm{M}_1\times \mathrm{M}_2$ imposing
that $\mathrm{K}_{\tau}$ is compatible with the transverse
holomorphic structure of $\mathcal{F}$ and that $\chi$ is of type
$(1,0)$. Using Newlander-Nirenberg theorem one can check that
$\mathrm{K}_\tau$ is integrable, i.e. a complex structure, if and
only $d\chi^{0,2}=0$, which holds in our case. Moreover the vector
field $v=T_1-\tau T_2$ is of type $(1,0)$, therefore it is
holomorphic if and only if $[v,Q^{0,1}]\subset Q^{0,1}$ where
$Q^{0,1}$ denotes the subbundle of $T(\mathrm{M}_1\times
\mathrm{M}_2)^{\C}$ of vector fields of type $(0,1)$ with respect to
$\mathrm{K}_{\tau}$, which holds as a consequence of the equality
$i_vd\chi=0$.
\end{proof}

\begin{exam}
Let $\mathrm{K}$ be a compact connected real Lie group of odd
dimension. Since $\mathrm{K}$ admits a nacs the previous proposition
describes a complex structure on the product $\mathrm{K}\times S^1$.
As $\mathrm{K}\times S^1$ is also a Lie group this can be seen as a
particular case of a result by Samelson (cf. \cite{Samel}) which
states that every compact connected Lie group of even dimension
admits a left-invariant complex structure. Nevertheless one obtains
more complex structures by this construction. Indeed, when
$\dim_{\R}\mathrm{(K)}>1$ we can assume that the nacs is
non-invariant. Moreover, by topological reasons $\mathrm{K}\times
S^1$ cannot be K\"{a}hlerian except when $\mathrm{K}$ is a real
torus.
\end{exam}

\begin{prop}\label{prop:caseB}\emph{\textbf{(Case B)}}
Let $\mathrm{M}$ be a manifold endowed with a nacs and ${\mathcal
F}$ its associated flow. Let $\pi:\mathrm{X}\rightarrow \mathrm{M}$
be a $S^1$-principal bundle over $\mathrm{M}$ with Chern class
$[d\beta]$, where $\beta$ is a 1-form on $\mathrm{X}$ such that
$d\beta$ is the pull-back of a closed $(1,1)$-form on $\mathrm{M}$
basic for the flow $\mathcal{F}$. Then there exists a 1-parametric
family of complex structures $K_{\tau}$ on $\mathrm{X}$ for $\tau\in
\C\backslash \R$ so that the complex manifold $\mathrm{X}$ admits a
non-vanishing holomorphic vector field $v$.
\end{prop}

\begin{proof}
Let $\omega$ be the characteristic 1-form of the nacs and let $T$ be
the vector field inducing the CR-action. We denote by
$\widetilde{T}$ the vector field on $\mathrm{X}$ contained in
$\ker\beta$ such that $\pi_{\ast}(\widetilde{T})=T$ and we define
the 1-form $\widetilde{\omega}=\pi^{\ast}\omega$. Let $R$ denote the
fundamental vector field of the action corresponding to the
$S^1$-principal bundle $\pi:\mathrm{X}\rightarrow \mathrm{M}$ such
that $\beta(R)=1$. We apply now the same arguments as in proposition
\ref{prop:caseA} using the vector fields $\widetilde{T}$ and $R$,
the distribution ${\mathcal D }=\ker\beta\cap
\ker\widetilde{\omega}$ and the transverse holomorphic structure for
the flow $\widetilde{\mathcal F}=\langle\widetilde{T},R\rangle$
induced by the CR-structure of $\mathrm{M}$.  We define a
complex-valued 1-form $\chi$ by imposing $\ker \chi=\mathcal{D}$,
$\chi(v)=1$ and $\chi(\overline{v})=0$. The holomorphic vector field
$v$ is $\widetilde{T}-\tau R$ for $\tau\in \C\backslash \R$. The
hypothesis $d\beta\in
\pi^{\ast}\Omega^{1,1}(\mathrm{M}/\mathcal{F})$ and $d\omega\in
\Omega^{1,1}(\mathrm{M}/\mathcal{F})$ imply that $d\chi$ is of type
$(1,1)$, thus the complex structure is integrable.
\end{proof}

\begin{defn}
Let $\mathrm{M}$ be a compact manifold with a nacs $(\Phi^{1,0},T)$.
We define
$$\mathrm{Aut}_{\mathcal{T}}(\mathrm{M})=\{f\in\mathrm{Aut}_{\mathrm{CR}}(\mathrm{M}):\,
f_{\ast}T=T\}.$$\index{$\mathrm{Aut}_{\mathcal{T}}(\mathrm{M})$}
\end{defn}

\begin{prop}\label{prop:caseC1}\emph{\textbf{(Case C)}}
Let $\mathrm{M}$ be a manifold endowed with a nacs
 $(\Phi^{1,0},T)$. Given $f\in\mathrm{Aut}_{\mathcal T}(\mathrm{M})$ the
suspension $\mathrm{X}$ of $\mathrm{M}$\index{suspension} by $f$, i.e.
$\mathrm{X}=\mathrm{M}\times \R/\sim$ where $(x,s)\sim(f(x),s+1)$, admits a
1-parametric family of complex structures $K_{\tau}$ for $\tau\in \C\backslash \R$
so that the complex manifold $\mathrm{X}$ admits a non-vanishing holomorphic vector
field $v$ induced by $T-\tau \derpp{s}$.
\end{prop}

\begin{proof}
The proof is straightforward using the same arguments as in the
previous cases, the distribution ${\mathcal D}$ is induced by the
CR-structure $\Phi^{1,0}$ on $\mathrm{M}$.
\end{proof}

\begin{defn}\label{def:complex}
Let $\mathrm{M}$ be a manifold endowed with a nacs, $T$ the vector field defining
the CR-action and $\mathcal{F}$ its associated flow. We say that the pair
$(\mathrm{X},v)$ of a compact complex manifold $\mathrm{X}$ and a non-singular
holomorphic vector field $v$ is a \emph{compact complexification} of the pair
$(\mathrm{M},T)$ if:
\begin{enumerate}[\bf (i)]
\item $\mathrm{M}$ is a real submanifold of $\mathrm{X}$.
\item The CR-structure of $\mathrm{M}$ is compatible with
the complex structure of $\mathrm{X}$.
\item There exists $\lambda\in \C$ such that $\RE(\lambda v)=T$.
\end{enumerate}
\end{defn}

Both the constructions of case A and case C produce complex
manifolds that are compact complexifications of the departing normal
almost contact manifolds. Indeed, if $v=T_1-\tau T_2$ for $\tau\in
\C \backslash \R$ then $T_1=\RE\big(\frac{\mathrm{i}\bar{\tau}}{\IM
\tau}v\big)$ and $T_2=\RE \big(\frac{\mathrm{i}}{\IM \tau}v\big)$.

\begin{thm}\label{teo:Rio1}
Every compact K\"{a}hler manifold admitting a non-vanishing holomorphic vector field
can be obtained as a suspension by of proposition \ref{prop:caseC1} (case C). In
particular it is a compact complexification of a compact manifold endowed with a
nacs.
\end{thm}

\begin{proof}
Let $X$ be a compact K\"{a}hler manifold admitting a non-vanishing
holomorphic vector field $v$. By a result by Carell-Lieberman (c.f.
\cite{CarLieb}) there exists a holomorphic 1-form $\chi$ over
$\mathrm{X}$ such that $\chi(v)\neq 0$. As $\mathrm{X}$ is compact
we can assume $\chi(v)=1$ and as $\mathrm{X}$ is K\"{a}hlerian we
have $d\chi=0$. Assume $b_1(\mathrm{X})=2k$. Let
$\gamma_1,...,\gamma_{2k}$ be closed paths giving a basis of
$H_1(\mathrm{X},\Z)$ modulus torsion and let $\xi_1,...,\xi_{2k}$ be
the dual basis of closed 1-forms. Fix a basis
$\omega_1,...,\omega_k$ of $H^0(\mathrm{X},\Omega^1)$. By Hodge's
decomposition theorem we have
$$\xi_i=a_i^1\omega_1+...+a_i^k\omega_k+
b_i^1\overline{\omega}_1+...+b_i^k\overline{\omega}_k+df_i=\eta_i+df_i$$
where $f_i$ is a differentiable function, for $i=1,...,2k$, and
$a_i^j,b_i^j\in \C$. By Stokes theorem the two sets of 1-forms $
\{\xi_i\}$ and $\{\eta_i\}$ have the same periods. In particular
$\{\eta_1,...,\eta_{2k}\}$ is a basis of $H^1(\mathrm{X}, \C)$ dual
of $\{\gamma_1,...,\gamma_{2k}\}$. Since $\Q+\mathrm{i}\Q$ is a
dense subset in $\C$ we can choose $a_i\in \C$ for $i=1,...,2k$
arbitrarily small so that $\eta=\chi+\sum a_i\eta_i$ is a closed
1-form and $\int_{\gamma_j}\eta\in \Q+\mathrm{i}\Q$ for
$j=1,...,2k$. Moreover by construction the 1-form $\eta$ is of the
form:
$$\eta=c_1\omega_1+...+c_k\omega_k+d_1\overline{\omega}_1+...+d_k\overline{\omega}_k$$
with $c_i,d_i\in \C$ for $i=1,...,k$. It follows that $\eta(v)$ is
constant and close to 1 by construction, set $\eta(v)=\delta$. In an
analogous way $\eta(\bar{v})$ is constant and close to 0, set
$\eta(\bar{v})=\epsilon$. Therefore $\Gamma=\big\{\int_{\gamma}\eta:
\gamma\in H_1(\mathrm{X},\Z)\big\}$ is finitely generated and it is
contained in $\Q+\mathrm{i}\Q$, thus $\Gamma\cong \Z+\mathrm{i}\Z$.
Fixing a base point $p_0$ the differentiable map
$$\pi_1:\mathrm{X}\rightarrow \C/\Gamma$$
$$\phantom{aaaahhhhh!} p\mapsto \int_{p_0}^p \eta \ \mathrm{mod} \ \Gamma $$
over the elliptic curve $\C/\Gamma$ is well defined. Furthermore
$\pi_1$ is a proper submersion and thus a fibration. The real vector
fields $v+\bar{v}$ and $\mathrm{i}(v-\bar{v})$ are transverse to the
fibres of $\pi_1$ and preserve the fibration (because
$\eta(v+\bar{v})=\epsilon+\delta$,
$\eta(\mathrm{i}(v-\bar{v}))=\mathrm{i}(\delta-\epsilon)$ are
constants close to 1 and $\mathrm{i}$ respectively and $\eta$ is
closed). Since $\eta(v+\bar{v})=\delta+\epsilon\sim 1$ we can find a
linear map $h:\C \rightarrow \R$ such that $h(\eta(v+\bar{v}))>0$
and $h(\Gamma)\subset \Z$. Let $\bar{h}:\C/\Gamma\rightarrow \R/\Z$
be the induced fibration. The composition $\pi_2=\bar{h}\circ
\pi_1:\mathrm{X}\rightarrow \R/\Z$ is a fibration over the circle.
The fibres of $\pi_2$, denoted by $\mathrm{M}=\pi_2^{-1}(p)$, admit
a CR-structure induced by the complex structure on $\mathrm{X}$.
There exists $a\in \C \backslash \R$ such that the real vector field
$v_1=\RE(av)$ is tangent to $\mathrm{M}$. As the flow associated to
$v$ is holomorphic the vector field $v_1$ preserves the CR-structure
of $\mathrm{M}$ and induces a transverse CR-action. On the other
hand there exists $b\in \R^{+}$ such that the vector field
$v_2=\RE(bv)$ projects over the vector field $\derpp{t}$ on $S^1$.
The flow of $v_2$ preserves the CR-structure over $\mathrm{M}$ and
clearly $[v_1,v_2]=0$. Finally setting $\tau=\bar{a}\cdot b^{-1}$ we
obtain $v_1-\tau v_2=\RE(av)-\tau \RE(bv)=\mu \cdot v$, where
$\mu\in \C$. Taking the automorphism $f$ over $\mathrm{M}$ induced
by the flow of $v_2$ for time $1$ the compact complexification of
case C for the preceding $\tau$ gives rise to the original complex
structure.
\end{proof}

We have shown the suspension of a complex manifold $\mathrm{N}$ by an automorphism
$g\in\mathrm{Aut}_{\C}(\mathrm{N})$ admits a nacs. Applying the construction of case
C to such a manifold is equivalent to consider the quotient $\mathrm{X}$ of
$\mathrm{N}\times \C$ by the subgroup generated by $F(x,z)=(f(x),z+1)$ and
$G(x,z)=(g(x),z+\tau)$, where $f\in \mathrm{Aut}_{\C}(\mathrm{N})$ so that $f\circ
g=g\circ f$, which we will call \emph{double suspension}. In this case there is a
holomorphic fibration $\pi:\mathrm{X}\rightarrow E_{\tau}=\C/\langle 1,\tau\rangle$
such that the vector field induced by $\derpp{z}$ is transverse to the fibers.

Assume that $\mathrm{X}$ is a compact K\"{a}hler manifold admitting
a holomorphic vector field $v$ without zeros and $\mathcal{F}$ the
flow induced by $v$. In the previous theorem we have seen that as a
consequence of Carrell-Lieberman theorem (c.f. \cite{CarLieb}) there
exists a closed holomorphic 1-form $\chi$ on $\mathrm{X}$ such that
$\chi(v)=1$. Then the complex structure on $\mathrm{X}$ is the only
one compatible with the transverse holomorphic structure of
$\mathcal{F}$ such that $\chi$ is of type $(1,0)$.

\begin{thm}\label{teo:Rio2} Let $\mathrm{X}$ be a
compact K\"{a}hler manifold admitting a holomorphic vector field $v$ without zeros
and $\mathcal{F}$ the flow induced by $v$ and let $\chi$ be a holomorphic 1-form
such that $\chi(v)=1$. Then there exists a closed 1-form $\beta$ of type (0,1)
arbitrarily small such that the complex structure on $\mathrm{X}$ compatible with
the transverse holomorphic structure of $\mathcal{F}$ and such that $\chi+\beta$ is
of type $(1,0)$ is induced by a double suspension. In particular every compact
K\"{a}hler manifold $\mathrm{X}$ admitting a holomorphic vector field $v$ without
zeros admits a complex structure induced by a double suspension on the underlying
smooth manifold $\mathrm{X}$ arbitrarily close to the original one.
\end{thm}

\begin{proof}
We proceed as in theorem \ref{teo:Rio1} to obtain from $\chi$ a
closed 1-form $\eta$ with group of periods $\Gamma\cong
\Z+\mathrm{i}\Z$ and a smooth fibration $\pi_1:\mathrm{M}\rightarrow
\C/\Gamma$ given by $x\mapsto \int_{x_0}^x \eta$. Every fibre
$\mathrm{N}$ of $\pi_1$ is transverse to the foliation
$\mathcal{F}_v$ generated by $v$. Therefore $\mathrm{N}$ admits a
complex structure. Note that $\eta=\pi_1^{\ast}(dz)$. Consider the
universal covering $p:\C \rightarrow \C/\Gamma_{\tau}$ and the
pullback $\pi_2:\mathrm{N}\times \C \rightarrow \C$ of the fibration
$\pi_1$ by the map $p$. There exists a map $q:\mathrm{N}\times \C
\rightarrow \mathrm{X}$ such that $\pi_1\circ q=p\circ \pi_2$. The
holomorphic vector field $v$ is transverse to the leaves of $\pi_1$
and it preserves the complex structure on $\mathrm{N}$. We recall
that $\eta(v)=\delta\sim 1$ and that $\eta(\bar{v})=\epsilon\sim 0$.
Fix $\tau \in \C \backslash \R$. We decompose $\delta^{-1}v=v_1-\tau
v_2$, where $v_1$ and $v_2$ are real vector fields. Then $v_1$ and
$v_2$ are transverse to the fibers of $\pi_1$, they preserve the
fibration and the complex structure on $\mathrm{N}$ and
$[v_1,v_2]=0$ (for $\delta^{-1}v$ is holomorphic). Finally, they
project over $\derpp{t}$ and $\derpp{s}$ on $\C/\Gamma$ respectively
where $z=\frac{\mathrm{i}}{2\IM \tau}(ds+\bar{\tau}dt)$. We set $f,g
\in \mathrm{Aut}_{\C}(\mathrm{N})$ as the flows $v_1$ and $v_2$ for
time $1$ and $-1$ respectively. Thus $\mathrm{X}$ is diffeomorphic
to the double suspension $\mathrm{N}\times \C/\langle F,G \rangle$
where $F(x,z)=(f(x),z+1)$, $G(x,z)=(g(x),z+\tau)$ and $f\circ
g=g\circ f$. Moreover the complex structure on $\mathrm{X}$ induced
by it which is arbitrarily close to the original one. Note that with
the new complex structure on $\mathrm{X}$ the fibration \linebreak
$\pi_1:\mathrm{X}\rightarrow \C/\Gamma$ is holomorphic and the
fibres $\mathrm{N}$ are analytic submanifolds. As we can choose
$\eta$ arbitrarily close to the starting holomorphic 1-form $\chi$
we can conclude.
\end{proof}

There is a natural generalization to a suspension of a compact complex manifold
$\mathrm{N}$ by a commutative subgroup $\Gamma=\langle f_1,...,f_s,g_1,...,g_s
\rangle$ of $\mathrm{Aut}_{\C}(\mathrm{N})$. The resulting manifold has dimension
$\dim_{\C}\mathrm{N}+s$ and fibers over the torus $\mathbb{T}^s$. Then one can
prove, with the same arguments as in theorem \ref{teo:Rio2}, the result below.

\begin{thm}\label{teo:Rio2+} Let $\mathrm{X}$ be a compact K\"{a}hler manifold such that
$\mathfrak{h}$ admits an abelian subalgebra
$\widetilde{\mathfrak{h}}$ of holomorphic vector fields without
zeros such that $\dim_{\C}\widetilde{\mathfrak{h}}=s>0$. The
underlying smooth manifold $\mathrm{X}$ admits a complex structure
arbitrarily close to the original one, obtained as a suspension over
the complex torus $\mathbb{T}^s$.
\end{thm}

We denote by $\mathfrak{h}$ the Lie algebra of holomorphic vector
fields on $\mathrm{X}$ and by $\mathfrak{h}_0$ the Lie algebra of
holomorphic vector field with zeros. Recall that
$[\mathfrak{h},\mathfrak{h}]\subset\mathfrak{h}_0$. Therefore if
$\dim_{\C}\mathfrak{h}=s>0$ and $\mathfrak{h}_0=0$ the hypothesis of
the above theorem hold. The limit case, i.e. when
$\dim_{\C}\widetilde{\mathfrak{h}}=\dim_{\C}\mathrm{X}$, is a
classical result by Wang's:

\begin{cor}\label{cor:wang}
Let $\mathrm{X}$ be a complex parallisable compact K\"{a}hler
manifold, then $\mathrm{X}$ is a complex torus.
\end{cor}

\subsection{Nacs on compact 3-manifolds}\label{sec:3nacs}

When M is a compact manifold of dimension 3 there is a classification due to
H.Geiges of the manifolds admitting a nacs based on the classification of compact
complex surfaces (see \cite{Gei}). Using the classification of transversely
holomorphic flows on a compact connected 3-manifold (see \cite{Bru} and
\cite{Ghys3}) together with the condition of the existence of a CR-structure and a
transverse CR-action we give an alternative proof of the classification. The main
interest of this point of view is that it determines not only the 3-manifold but
also the flow of the CR-action.

\begin{prop}
Let $\mathrm{M}$ be a compact connected 3-manifold endowed with a nacs. Then, up to
diffeomorphism, the manifold $\mathrm{M}$ and the vector field inducing the
\emph{CR}-action belong to the following list:
\begin{enumerate}[\bf i)]
\item Seifert fibrations over a Riemann surface such that the isometric flow of the
$S^1$-action admits a characteristic 1-form $\omega$ such that $d\omega$ is of type
$(1,1)$. \item Linear vector fields in $\mathbb{T}^3$. \item Flows on $S^3$ induced
by a singularity of a holomorphic vector field in $\C^2$ in the Poincar\'{e} domain
and their finite quotients, i.e. flows on lens spaces. \item Suspensions of a
holomorphic automorphism of $\pp^1$ with a vector field tangent to the flow
associated to the suspension.
\end{enumerate}
Moreover, all the previous manifolds admit a nacs such that the \emph{CR}-action is
the one induced by the corresponding vector field.
\end{prop}

\begin{proof}
Since the last statement of the proposition is clear it is enough to
depart from the classification of transversely holomorphic flows and
rule out the two cases that do not admit a nacs: strong stable
foliations associated to suspensions of hyperbolic diffeomorphisms
of $\mathbb{T}^2$ and $\C\times \R\backslash \{(0,0)\}/\sim$ where
$(z,t)\sim (\lambda z,2t)$ for $\lambda \in \C$ such that
$|\lambda|>1$ with the flow induced by the vertical vector field
$\derpp{t}$. The first ones are are examples of non-isometric
Riemannian flows (see \cite{Carr}), if they admitted an invariant
CR-structure together with a transverse CR-action they would be
isometric (for the flow of a nacs admits a transverse invariant
distribution). Let now $\mathrm{M}$ be $\C\times \R\backslash
\{(0,0)\}/\sim$ where $(z,t)\sim (\lambda z,2t)$ for $\lambda \in
\C$ such that $|\lambda|>1$ with the flow $\mathcal{F}$ induced by
the vertical vector field $\derpp{t}$. Suppose that $\mathrm{M}$
admits a CR-structure $\Phi^{1,0}$ transverse to ${\mathcal F}$ and
a vector field $T$ tangent to ${\mathcal F}$ inducing a CR-action.
Then $\mathrm{M}^3\times S^1$ is a compact complex surface that
admits a holomorphic non-vanishing vector field. Therefore it
belongs to the following list (c.f. \cite{DlOeTo}, \cite{DlOeTo2}):
\begin{enumerate}[\bf I)]
\item Complex tori.
\item Principal Seifert fibre bundles over a Riemann surface of genus $g\geq 1$ with fiber
an elliptic curve.
\item Ruled surfaces over an elliptic curve.
\item Almost-homogeneous Hopf-surfaces.
\end{enumerate}
Since $\mathrm{X}$ is homeomorphic to $S^2\times S^1\times S^1$ the
complex surface $\mathrm{X}$ must be a ruled surface over an
elliptic curve. However we will see that this is a contradiction.
Recall that the universal covering of a ruled surface over an
elliptic curve is either $\D \times \pp^1$ or $\C \times \pp^1$. By
construction of the complex structure on
$\mathrm{X}=\mathrm{M}\times S^1$ we have that $t+\tau s$ is a
holomorphic coordinate for some $\tau\in \C\backslash \R$, therefore
the analytic universal covering $\widetilde{\mathrm{X}}=\C \times \R
\backslash \{(0,0)\}\times \R$ of $\mathrm{X}$ admits a holomorphic
projection $p:\widetilde{\mathrm{X}}\rightarrow \C$ defined by
$p(z,t,s)=z$ with fiber an open subset of $\C$. As $\pp^1$ is
compact by the maximum principle it is immersed in the fibers of
$p$, which is a contradiction.
\end{proof}

\section{K\"{a}hler criteria}

\subsection{Conditions on the Euler class}

\begin{thm}\label{thm:euler}
Let $\mathrm{X}$ be a compact complex manifold with a non-vanishing
holomorphic vector field $v$. For every $\tau\in \C\backslash \R$
let $T_1$ and $T_2$ be the two real vector fields defined by
$v=T_1-\tau T_2$ and $\mathcal{F}_1$, $\mathcal{F}_2$ the flows
defined by $T_1$, $T_2$ respectively. If $\mathrm{X}$ is
K\"{a}hlerian then
$e_{\mathcal{F}_1}(\mathrm{X})=e_{\mathcal{F}_2}(\mathrm{X})=0$.
\end{thm}

\begin{proof}
Since $\mathrm{X}$ is a compact K{\"a}hler manifold with a
holomorphic vector field $v$ without zeros by Carrell-Liebermann's
theorem there exists a closed holomorphic 1-form $\alpha$ such that
$\alpha(v)=1$. We decompose
$\alpha=\frac{i}{2\IM\tau}(\alpha_2+\bar{\tau}\alpha_1)$ where
$\alpha_1$ and $\alpha_2$ are real closed 1-forms. Using that
$\alpha(v)=1$ and $\alpha(\overline{v})=0$, a direct computation
shows that $\alpha_i(T_j)=\delta_{ij}$ for $i,j=1,2$. Thus
$e_{\mathcal{F}_1}(\mathrm{X})=e_{\mathcal{F}_2}(\mathrm{X})=0$.
\end{proof}

\begin{prop}\label{prop:complex}
Let $\mathrm{M}$ be a compact manifold endowed with a nacs, let $T$
be the vector field of the \emph{CR}-action and $\mathcal{F}$ its
associated flow. Assume that the compact complex manifold
$\mathrm{X}$ together with a holomorphic vector field $v$ are a
compact complexification of $(\mathrm{M},T)$. If $\,\mathrm{X}$ is
K\"{a}hlerian then $e_{\mathcal{F}}(\mathrm{M})=0$.
\end{prop}

\begin{proof}
We denote by $v$ the holomorphic vector field on $\mathrm{X}$ and by
$\lambda$ the complex number such that $T=\RE(\lambda v)$ on
$\mathrm{M}$. Since $\mathrm{X}$ is a compact K{\"a}hler manifold
with a non-singular vector field $v$ there exists a holomorphic
closed 1-form $\alpha$ such that $\alpha(v)=\lambda^{-1}$. We can
decompose $(\lambda\cdot v)_{|\mathrm{M}}=T-iS$ where $S$ is a real
vector field. Set $\alpha=\frac{1}{2}(\alpha_1+i\alpha_2)$ where
$\alpha_1,\alpha_2$ are real 1-forms on $\mathrm{X}$. Then
$\alpha_1$ and $\alpha_2$ are closed and $\alpha_1(T)=1$ on
$\mathrm{M}$. The closed real 1-form
$\omega:={\alpha_1}_{|\mathrm{M}}$ verifies $\omega(T)=1$ (because
$\alpha(v)=1$ and $\alpha(\overline{v})=0$), therefore
$e_{\mathcal{F}}(\mathrm{M})=0$.
\end{proof}

\begin{cor}
Let $\mathrm{M}$ be a compact manifold endowed with a nacs. If
$e_{\mathcal{F}}(\mathrm{M})\neq 0$ then no compact complexification
obtained as a product by proposition \ref{prop:caseA} (case A) or as
a suspension by proposition \ref{prop:caseC1} (case C) is
K\"{a}hlerian.
\end{cor}

\begin{cor}
Let $\mathrm{M}$ be a compact manifold endowed with a nacs, let $T$
be the vector field of the \emph{CR}-action and $\mathcal{F}$ its
associated flow. If $\,\mathrm{M}$ admits a normal contact structure
compatible with the \emph{CR}-action induced by $T$ then
$(\mathrm{M},T)$ admits no K{\"a}hler compact complexification.
\end{cor}

It follows from corollary \ref{cor:contacte}.

\begin{cor}
Let $\mathrm{M}$ be a compact manifold endowed with a nacs. If
$b_1(\mathrm{M})=0$, in particular if $\,\mathrm{M}$ is simply
connected, then $(\mathrm{M},T)$ admits no K{\"a}hler compact
complexification\index{compact complexification}.
\end{cor}

It is a consequence of proposition \ref{cor:tischler}.

\begin{cor}
Let $\,\mathrm{M}$ be a compact connected semi\-simple real Lie group of odd
dimension endowed with a nacs and $T$ the vector field of the \emph{CR}-action. Then
$(\mathrm{M},T)$ admits no K{\"a}hler compact complexification.
\end{cor}

It is enough to recall that for any such group $b_1(\mathrm{M})=0$.

\begin{prop}\label{prop:caseBEulerKahler}
Let $\mathrm{M}$ be a compact manifold endowed with a nacs and let $\mathrm{X}$ be a
compact complex manifold obtained as a $S^1$-principal bundle
$\pi:\mathrm{X}\rightarrow \mathrm{M}$ over $\mathrm{M}$ by proposition
$\ref{prop:caseB}$ (case B). If $\,\mathrm{X}$ is K\"{a}hlerian then
$e_{\mathcal{F}}(\mathrm{M})=0$ and the $S^1$-principal bundle is flat. In
particular, if $\mathrm{X}$ is K\"{a}hlerian and $H^2(\mathrm{M},\Z)$ has no torsion
then the $S^1$-principal bundle is topologically trivial. Moreover, if $\alpha$ is a
connection 1-form on $\mathrm{X}$ such that $d\alpha\in
\pi^{\ast}\Omega^{1,1}(\mathrm{M}/\mathcal{F})$ then $[d\alpha]=0$ in
$H^2(\mathrm{M}/\mathcal{F})$.
\end{prop}

\begin{proof}
If $v$ is the holomorphic vector field of the complexification there
exists a closed holomorphic 1-form $\alpha$ on $\,\mathrm{X}$ such
that $\alpha(v)=1$. The connected group $S^1$ acts holomorphically
on $\mathrm{X}$ (as the group of the action of the $S^1$-principal
bundle), therefore the forms $\alpha$ and $\bar{\alpha}$ are
invariant by the action of $S^1$. Notice that $v=\widetilde{T}-\tau
R$ where $R$ is the vector field of the action and, $\widetilde{T}$
is the vector field contained in $\ker \beta$ such that
$\pi_{\ast}(\widetilde{T})=T$ and $\tau\in \C\backslash \R$. We
decompose $\alpha=\frac{i}{2\IM \tau}(\alpha_2+\bar{\tau}\alpha_1)$
where $\alpha_1,\alpha_2$ are real 1-forms. Then $\alpha_1$ and
$\alpha_2$ are closed 1-forms invariant by the action of $S^1$ (for
they are a linear combination of $\alpha$, $\bar{\alpha}$ with
constant coefficients) such that
$\alpha_1(\widetilde{T})=\alpha_2(R)=1$ and
$\alpha_1(R)=\alpha_2(\widetilde{T})=0$ (because $\alpha(v)=1$ and
$\alpha(\overline{v})=0$). Since $\alpha_1$ is a closed real basic
$S^1$-invariant 1-form it induces a closed 1-form $\omega$ on
$\mathrm{M}$ such that $\omega(T)=1$, thus
$e_{\mathcal{F}}(\mathrm{M})=0$. Finally, $\alpha_2$ is a closed
connection 1-form for the $S^1$-principal bundle
$\pi:\mathrm{X}\rightarrow \mathrm{M}$, so it is flat. When
$H^2(\mathrm{M},\Z)$ has no torsion flat bundles are topologically
trivial. Moreover, if $\alpha$ is a connection 1-form on
$\mathrm{X}$ such that $d\alpha\in
\pi^{\ast}\Omega^{1,1}(\mathrm{M}/\mathcal{F})$ then
$[d\alpha]=[d\alpha_2]=0$ in $H^2(\mathrm{M}/\mathcal{F})$.
\end{proof}

Finally for the suspension (case C) of a compact connected Lie group
endowed with a nacs an elementary computation of the second
cohomology group shows that the resulting complex manifold cannot be
K\"{a}hlerian:

\begin{prop}\label{prop:CLieGroups}
Let $\mathrm{K}$ be a non-abelian compact connected real Lie group of odd dimension
endowed with a nacs, $f\in \mathrm{Aut}_{\mathcal{T}}(\mathrm{K})$ and
$\mathrm{X}=\mathrm{K}\times_f \R$ endowed with the complex structure obtained as a
suspension by proposition \ref{prop:caseC1} (case C). Then $\mathrm{X}$ is not
K\"{a}hlerian.
\end{prop}

\begin{proof}
The suspension $\mathrm{X}=\mathrm{K}\times_f \R$ admits a finite
covering $\widetilde{\mathrm{X}}=\mathrm{M}\times_{\tilde{f}} \R$
such that $\mathrm{M}\cong \mathrm{K}'\times (S^1)^r$ where
$\mathrm{K}'$ is a compact connected semisimple real Lie group,
$0\leq r<\dim_{\R} \mathrm{K}$ and $\widetilde{f}$ the lift of $f$
to $\widetilde{\mathrm{X}}$. Since
$b_1(\mathrm{K}')=b_2(\mathrm{K}')=0$ we conclude that
$H^2(\mathrm{M})\cong H^2((S^1)^r)$ and $H^1(\mathrm{M})\cong
H^1((S^1)^r)$. Then using Mayer-Vietoris sequence for the De Rham
cohomology groups (c.f. \cite{BottTu}) one proves that
$H^2(\widetilde{\mathrm{X}})$ is isomorphic to
$$\{[\sigma]\in H^2(\mathrm{M}):\, f^{\ast}[\sigma]=[\sigma] \}\oplus
\bigg(\frac{H^1(\mathrm{M})}{\{[\sigma-f^{\ast}\sigma]:\,
[\sigma]\in H^1(\mathrm{M})\}}\bigg)\wedge [ds].$$ It is not
difficult to see that $\widetilde{\mathrm{X}}$ can not be
K\"{a}hlerian and it follows that $\mathrm{X}$ is not K\"{a}hlerian.
\end{proof}

\pagebreak
\subsection{Criteria for isometric flows}

\begin{thm}\label{thm:kahleriso}
Let $\mathrm{X}$ be a compact complex manifold with a non-vanishing
holomorphic vector field $v$. For every $\tau\in \C\backslash \R$
let $T_1$ and $T_2$ be the two real vector fields defined by
$v=T_1-\tau T_2$ and $\mathcal{F}_1$, $\mathcal{F}_2$ the flows
defined by $T_1$, $T_2$ respectively. Assume that the flows
$\mathcal{F}_1$ and $\mathcal{F}_2$ are isometric. Then the manifold
$\mathrm{X}$ is K\"{a}hlerian if and only if $e_{{\mathcal
F}_1}(\mathrm{X})=e_{{\mathcal F}_2}(\mathrm{X})=0$ and the real
foliation ${\mathcal F}=\langle T_1,T_2\rangle$ is transversely
K\"{a}hlerian.
\end{thm}

\begin{proof}
$\Rightarrow):$ The same argument as in theorem \ref{thm:euler}
shows that there are two closed real 1-forms $\alpha_1$ and
$\alpha_2$ on $\mathrm{M}$ such that $\alpha_i(T_j)=\delta_{ij}$ for
$i,j=1,2$. Then $\alpha_1$ and $\alpha_2$ are characteristic forms
for ${\mathcal F}_1$ and ${\mathcal F}_2$ respectively and
$e_{{\mathcal F}_1}(\mathrm{X})=e_{{\mathcal F}_2}(\mathrm{X})=0$.
We denote by $(\varphi_1)_t$ and $(\varphi_2)_t$ the 1-parametric
groups associated to $T_1$ and $T_2$ respectively and by
$\mathrm{H}$ the closure in $\mathrm{Isom(X)}$ of the abelian group
generated by $(\varphi_1)_t$ and $(\varphi_2)_t$.
If $\Phi$ is the K\"{a}hler form on $\mathrm{X}$ then the transverse
part $\Psi(\cdot,\cdot)$ with respect to $\mathcal{F}$ of
$$\int_{\mathrm{H}}\Phi(\sigma_{\ast}\cdot,\sigma_{\ast}\cdot),$$
where we integrate with respect to the Haar mesure on $\mathrm{H}$,
is a transverse K\"{a}hler form.

$\Leftarrow):$ As $\mathcal{F}_1$ and $\mathcal{F}_2$ are isometric
we can assume that there exists an invariant transverse distribution
$\mathcal{D}$ of maximal dimension for the real foliation
$\mathcal{F}=\langle T_1,T_2 \rangle$. We denote by $\omega_1$ and
$\omega_2$ the 1-forms on $\mathrm{X}$ defined by
$\omega_i(T_j)=\delta_{ij}$ and ${\omega_i}_{|\mathcal{D}}=0$ for
$i,j=1,2$. Since $e_{{\mathcal F}_1}(\mathrm{X})=e_{{\mathcal
F}_2}(\mathrm{X})=0$ there exist $\beta_1,\beta_2\in
\Omega^1(\mathrm{X}/{\mathcal F})$ such that $d\beta_i=d\omega_i$
for $i=1,2$. We denote by $\beta$ the basic form
$\beta=\frac{i}{2\IM\tau}(\beta_2+\overline{\tau} \beta_1)$. It
follows that $d\beta=d\chi$. We begin by showing that it is enough
to find $\alpha\in \Omega^1(\mathrm{X}/{\mathcal F},\C)$ of type
$(1,0)$ such that $d\alpha=d\chi$. Indeed, if $\alpha$ exists the
form $\Phi=(\chi-\alpha)\wedge (\overline{\chi}-\overline{\alpha})$
is closed and of type $(1,1)$. Adding to $\Phi$ a positive multiple
of the transverse K\"{a}hler form of ${\mathcal F}$ we obtain a
K\"{a}hler form on $\mathrm{X}$ and the proof is complete. We will
now show that such a form exists. Since $d\beta^{0,2}=d\chi^{0,2}=0$
we have $d\beta=d(\beta^{1,0})+\partial (\beta^{0,1})$, i.e.
$\overline{\partial}(\beta^{0,1})=0$, and
$$d(\partial\beta^{0,1})=(\partial+\overline{\partial})(\partial
\beta^{0,1})= -\partial\overline{\partial}\beta^{0,1}=0$$ so
$\partial(\beta^{0,1})$ is a $(1,1)$-form which is $\partial$-exact
as a basic form and $d$-closed. Applying the basic
$\partial\overline{\partial}$-lemma (c.f. \cite{Kacimi}) to
$\partial(\beta^{0,1})$ we obtain a basic function $f$ such that
$$\partial(\beta^{0,1})=\partial\overline{\partial}f=\overline{\partial}\partial(-f).$$
Then $-\partial f$ is a basic form of type $(1,0)$ such that
$d(-\partial f)=-\overline{\partial}\partial
f=\partial(\beta^{0,1})$. The form $\alpha=\beta^{1,0}-\partial f$
is basic, of type $(1,0)$ and $d\alpha=d\beta=d\chi$ so the
conclusion follows.
\end{proof}

Every Riemannian holomorphic flow ${\mathcal F}$ in a compact complex surface $S$ is
transversely K\"{a}hlerian. Therefore with the above hypothesis when
$\dim_{\C}\mathrm{X}=2$ the complex manifold $\mathrm{X}$ is K\"{a}hlerian if and
only if $e_{{\mathcal F}_1}(\mathrm{X})=e_{{\mathcal F}_2}(\mathrm{X})=0.$

\begin{cor}
Let $\mathrm{X}$ be a complex manifold obtained as a product by proposition
\ref{prop:caseA} (case A) from two manifolds $\mathrm{M}_1$ and $\mathrm{M}_2$
endowed with a nacs such that the flows $\mathcal{F}_1$ and $\mathcal{F}_2$ in
$\mathrm{M}_1$ and $\mathrm{M}_2$ respectively associated to the nacs are isometric.
Then $\mathrm{X}$ is K\"{a}hlerian if and only if
$e_{\mathcal{F}_1}(\mathrm{M}_1)=e_{\mathcal{F}_2}(\mathrm{M}_2)=0$ and the flows
$\mathcal{F}_1$ and $\mathcal{F}_2$ are transversely K{\"a}hlerian.
\end{cor}

\begin{cor}
Let $\mathrm{M}$ be a compact manifold endowed with a nacs such that its associated
flow $\mathcal{F}$ is isometric. Let $\mathrm{X}$ be a complex manifold obtained as
a $S^1$-principal bundle over $\mathrm{M}$ by proposition \ref{prop:caseB} (case B).
Then $\mathrm{X}$ is K\"{a}hler if and only if the $S^1$-principal bundle
$\pi:\mathrm{X}\rightarrow \mathrm{M}$ is flat, the flow $\mathcal{F}$ is
transversely K{\"a}hlerian on $\mathrm{M}$ and $e_{\mathcal{F}}(\mathrm{M})=0$.
\end{cor}

\begin{exam}
Let $\mathrm{S}$ be a compact complex surface obtained as a $S^1$-principal bundle
over a 3-manifold $\mathrm{M}$ with a nacs by means of the construction of case B.
From the last result together with the classification of compact complex surfaces
one concludes the following. With the notation of section \S\ref{sec:3nacs} for each
of the possibilities for $\mathrm{M}$ the corresponding surface $\mathrm{S}$ is:
\textbf{i)} an elliptic Seifert principal fibre bundle, \textbf{ii)} a complex torus
or non-K\"{a}hlerian elliptic Seifert principal fibre bundle, \textbf{iii)} a Hopf
surface and \textbf{iv)} either a Hopf surface or a ruled surface over an elliptic
curve.
\end{exam}

\begin{prop}\label{prop:kahlersusp2}
Let $\mathrm{X}$ be a complex manifold obtained as a suspension by proposition
\ref{prop:caseC1} (case C) from a manifold $\mathrm{M}$ endowed with a nacs such
that its associated flow $\mathcal{F}$ is isometric. If $\mathrm{X}$ is
K\"{a}hlerian then the flow $\mathcal{F}$ is transverselly K\"{a}hlerian and there
exists a basic K\"{a}hler form $\Phi$ such that $[f^{\ast}\Phi]=[\Phi]\in
H^{1,1}(\mathrm{M}/{\mathcal F})$.
\end{prop}

\begin{proof}
We choose a K\"{a}hler form $\Psi_0$ on $\mathrm{X}$ and define
$\Psi$ as the pullback to $\mathrm{M}\times \R$. We define
$\Phi_0=\Psi(x,s_0)_{|\mathrm{M}}$ so that
$[f_{\ast}\Phi_0]=[\Phi_0]$ and $\Phi_0$ is a closed form of type
$(1,1)$ with respect to the holomorphic transverse structure. We
denote by $\varphi_t$ the 1-parameter group on $\mathrm{M}$
associated to $T$ and by $\mathrm{H}$ the closure in
$\mathrm{Isom(M)}$ of the abelian group generated by $\varphi_t$.
The form $\Phi$ defined as the basic part of
$$\int_{\mathrm{H}}\Phi_0(\sigma_{\ast}\cdot,\sigma_{\ast}\cdot),$$
where we integrate with respect to the Haar mesure on $\mathrm{H}$,
is a basic K\"{a}hler form on $\mathrm{M}$. Moreover since
$f^{\ast}\Phi_0=\Phi_0+d\alpha$ where $\alpha$ is a 1-form on
$\mathrm{M}$, the form $f^{\ast}\Phi$ is the basic part of
\begin{eqnarray*}
\int_{\mathrm{H}} f^{\ast}(\Phi_0 (\sigma_{\ast}
\cdot,\sigma_{\ast}\cdot))&=&\int_{\mathrm{H}}(f^{\ast}\Phi_0)(\sigma_{\ast}\cdot,\sigma_{\ast}\cdot)=
\int_{\mathrm{H}}\Phi_0(\sigma_{\ast}\cdot,\sigma_{\ast}\cdot)+
\int_\mathrm{H} d\alpha(\sigma_{\ast}\cdot,\sigma_{\ast}\cdot)
\\&=&\int_{\mathrm{H}}\Phi_0(\sigma_{\ast}\cdot,\sigma_{\ast}\cdot)+ d \int_\mathrm{H} \alpha(\sigma_{\ast}\cdot)
\end{eqnarray*}
(see \cite{Greub1} and \cite{Greub2} and note that $f_{\ast}T=T$ so
$f\circ \sigma=\sigma \circ f$). Therefore $[f^{\ast}\Phi]=[\Phi]$.
\end{proof}

\subsection{Suspensions of Levi-flat nacs}

In this section we find necessary and suficient conditions for a complex manifold
obtained as a suspension by the construction of case C from a manifold endowed with
a Levi-flat nacs to be K\"{a}hlerian.

Let $\mathrm{M}$ be a compact manifold endowed with a nacs and let
$\omega$ be its characteristic 1-form. Assume that the CR-structure
is Levi-flat, i.e. $d\omega=0$. The suspension $\mathrm{M}$ of a
compact complex manifold $\mathrm{N}$ by $g\in Aut_{\C}(\mathrm{N})$
with the natural nacs is an example of this situation. Conversely,
one has the following result:

\begin{prop}
Let $\mathrm{M}$ be a manifold endowed with a nacs with Levi-flat
\emph{CR}-structure. If its characteristic 1-form $\omega$ has group of periods
$\Gamma_{\omega}\cong \Z$ then $\mathrm{M}$ and the given nacs can be obtained as a
suspension of a compact complex manifold $\mathrm{V}$ by $g\in
\mathrm{Aut}_{\C}(\mathrm{V})$.
\end{prop}

\begin{proof} Since $\omega$ has periods group $\Gamma_{\omega}\cong
\Z$ there is a well defined fibration
$$\pi:\mathrm{M}\rightarrow S^1\cong \R/\Gamma_{\omega}, \quad x\mapsto \int_{x_0}^x \omega$$
with compact fibers isomorphic to $\mathrm{V}$ which are the leaves
of the foliation induced by $\ker\omega$. The compact leaf
$\mathrm{V}$ carries a complex structure induced by the CR-structure
$\Phi^{1,0}$ on $\mathrm{M}$. The automorphism $g\in
\mathrm{Aut}(\mathrm{V})$ corresponding to the suspension giving
rise to $\mathrm{M}$ is induced by $\varphi_{t_0}$, where $t_0$
verifies $\varphi_t(x_0)\notin V$ for $0<t<t_0$ and
$\varphi_{t_0}(x_0)\in V$ for every $x_0\in V$. The vector field of
the CR-action $T$ is transverse to the leaves of the fibration and
$\omega=\pi^{\ast}(dt)$, where $t$ denotes the real coordinate on
$S^1$. The choice of $g$ and the hypothesis $\omega(T)=1$ imply that
$T$ is the vector field induced by $\derpp{t}$.
\end{proof}

\begin{thm}\label{thm:kahlersusp2}
Let $\mathrm{X}$ be a complex manifold obtained as a suspension by
proposition \ref{prop:caseC1} (case C) from a manifold $\mathrm{M}$
endowed with a nacs with Levi-flat \emph{CR}-structure and $f\in
\mathrm{Aut}_{\mathcal{T}}(\mathrm{M})$. Assume that its associated
flow $\mathcal{F}$ is isometric and that $f^{\ast}=\mathrm{id}$
acting on $H^1(\mathrm{M},\C)$. Then $\mathrm{X}$ is K\"{a}hlerian
if and only there exists a basic K\"{a}hler form $\Phi$ for the flow
$\mathcal{F}$ such that $\,[f^{\ast}\Phi]=[\Phi]\in
H^{1,1}(\mathrm{M}/{\mathcal F})$.
\end{thm}

\begin{lem}\label{lem:Blanch3}
Let $\mathrm{M}$ a compact manifold with a transversely
K\"{a}hlerian flow $\mathcal{F}$ depending smoothly on a real
parameter $s\in \R$. Let $\omega$ be an exact basic form on
$\mathrm{M}$ of type $(p,q)$ such that $\omega=d_\mathrm{M}\alpha$
for a basic form $\alpha$ depending smoothly on $s$. Then there
exists a basic $(p-1,q)$-form $\mu$ and a basic $(p,q-1)$-form $\nu$
that depend smoothly on $s$ and such that
$\omega=d_\mathrm{M}\mu=d_\mathrm{M}\nu$.
\end{lem}

\begin{lem}\label{lem:Blanch4}
Let $\mathrm{M}$ a compact manifold with a transversely
K\"{a}hlerian flow $\mathcal{F}$ depending smoothly on a real
parameter $s\in \R$. Let $\omega$ be an exact basic form on
$\mathrm{M}$ such that $\omega=d_\mathrm{M}\alpha$ for a basic form
$\alpha$ depending smoothly on $s$. Assume that
$\omega=\omega_1+\omega_2$ where $\omega_1$ is of type $(p+1,q)$ and
$\omega_2$ is of type $(p,q+1)$, both basic and depending smoothly
on $s$. There exists a basic $(p,q)$-form $\nu$ such that
$\omega=d_\mathrm{M}\nu$ and it depends smoothly on $s$.
\end{lem}

The proofs of the previous lemmas are analogous to the ones in
\cite{Blanch}, pages 185-186, using the fact that on a transversely
K\"{a}hler isometric flow there exists a transversal Hodge theory
regarding basic forms and in particular a basic
$\partial\overline{\partial}$-lemma (c.f. \cite{Kacimi}). Note that
locally every suspension is a product, thus we can consider the
exterior derivative $d_s$ with respect to the local coordinate $s$
and then $d=d_{\mathrm{M}}+d_s$.

\begin{proof}
\noindent $\Longrightarrow):$ It follows from proposition
\ref{prop:kahlersusp2}.

\noindent $\Longleftarrow):$ We define $\alpha$ as the (1,0)-form
induced by $\frac{\mathrm{i}}{2\IM \tau}(ds+\bar{\tau} \omega)$,
which is a closed holomorphic 1-form on $\mathrm{X}$ such that
$\alpha(v)=1$. Besides, from Gysin's exact sequence we conclude that
$$H^1(\mathrm{M})\cong H^1(\mathrm{M}/\mathcal{F})\oplus \langle [\omega] \rangle.$$
From the hypothesis $f^{\ast}=\mathrm{id}$ acting on
$H^1(\mathrm{X})$ and a short computation using Mayer-Vietoris
sequence for the De Rham cohomology groups (c.f. \cite{BottTu}), it
follows that $H^1(\mathrm{X})\cong H^1(\mathrm{M})\oplus \langle
[ds]\rangle$. Since the flow $\mathcal{F}$ is transversely
K\"{a}hlerian on $\mathrm{M}$ by Hodge theory
$$H^1(\mathrm{M}/\mathcal{F})\cong H^{1,0}(\mathrm{M}/\mathcal{F})\oplus H^{0,1}(\mathrm{M}/\mathcal{F})$$
and there exists a basis $\alpha_1,...,\alpha_k$ of
$H^{1,0}(\mathrm{M}/{\mathcal{F}})$ of closed (1,0)-forms (and
therefore holomorphic) on $\mathrm{M}$ such that
$f^{\ast}\alpha_j=\alpha_j$ for $j=1,...,k$ (for
$f^{\ast}=\mathrm{id}$ on $H^1(\mathrm{M})$ and the
repre\-sen\-tatives of cohomology classes of type (1,0) are unique
as a consequence of the transversal Hodge theory). Note also that
$H^{1,0}(\mathrm{X}/\mathcal{F}_v)\cong
H^{1,0}(\mathrm{M}/\mathcal{F})$ since a basic 1-form on
$\mathrm{X}$ can not depend on $s$ and therefore lives on
$\mathrm{M}$. Fix the transverse K\"{a}hler form $\Phi$ on
$\mathrm{M}$. We can choose an open covering $\{U_i\}$ of $S^1$ so
that the fibration $p:\mathrm{X}\rightarrow S^1$ is trivial over
$U_i$. Then $\Phi$ induces a well-defined form $\Phi_i$ on
$p^{-1}(U_i)\cong U_i\times \mathrm{M}$. Let $\{\rho_i\}$ be a unit
partition associated to $\{U_i\}$. Then $\Phi_0=\sum_i
\rho_i(s)\Phi_i$ is a real global $(1,1)$-form on $\mathrm{X}$ so
that ${\Phi_0}_{|\mathrm{M}}$ is a transverse K\"{a}hler form on
$\mathrm{M}$ representing a fixed cohomology class, since
$[f^{\ast}\Phi]=[\Phi]$. Moreover $d\Phi_0=d_s \Phi_0$. Now we want
to obtain a closed real valued $(1,1)$-form $\widetilde{\Phi}$ on
$\mathrm{X}$ such that $\widetilde{\Phi}_{|\mathrm{M}}=\Phi_{0}$. We
search $\widetilde{\Phi}$ of the type:
$$\widetilde{\Phi}=\Phi_0+ H\wedge \bar{\alpha}+ \overline{H}\wedge \alpha+\mathrm{i}F \alpha\wedge \bar{\alpha},$$
where $H$ is a basic $(1,0)$-form on $\mathrm{X}$ and $F$ a
real-valued function on $\mathrm{X}$. The hypothesis
$d\widetilde{\Phi}=0$ is then equivalent to the following equations:
$$\left\{\begin{array}{rcll}
d_s \Phi_0-\frac{1}{2\mathrm{i}\IM \tau}(-d_\mathrm{M}H\wedge ds+d_\mathrm{M} \overline{H}\wedge ds) &=&0  &\quad \mathrm{(I)}\\
-\tau d_\mathrm{M}H\wedge \omega+\bar{\tau} d_\mathrm{M} \overline{H}\wedge \omega&= & 0 & \quad \mathrm{(II)} \\
\frac{1}{2\mathrm{i}\IM \tau}(-\tau d_sH\wedge \omega+\bar{\tau} d_s
\overline{H}\wedge \omega)-idF \alpha\wedge \bar{\alpha}&=& 0& \quad
\mathrm{(III).}
\end{array}\right.$$
Note that $$\mathrm{i}dF \alpha\wedge \bar{\alpha}=\frac{-dF}{2\IM
\tau} ds\wedge \omega.$$ Obtaining the form $\widetilde{\Phi}$ is
equivalent to finding $H$ and $F$ that solve the previous system.
Roughly speaking, we will first solve (I) and (II) so that we
determine $H$ and then define $F$ as the solution of (III). The
equation (I) is equivalent to
$$2\mathrm{i}\IM \tau i_{\derpp{s}}d_s \Phi_0=d_\mathrm{M} \overline{H}-d_\mathrm{M} H=d_\mathrm{M}(\overline{H}-H)$$
which is fulfilled if
$$d_\mathrm{M}H=\bar{\tau} i_{\derpp{s}}d_s \Phi_0 \qquad \mathrm{(IV).}$$
The form $\gamma:=\bar{\tau} i_{\derpp{s}}d_s \Phi_0$ can be seen
locally as a basic form on $\mathrm{M}$ of type (1,1) that depend
smoothly on the real parameter $s$. To apply lemma \ref{lem:Blanch3}
to obtain local basic $(1,0)$-forms $H_i$ on $p^{-1}(U_i)\cong
U_i\times \mathrm{M}$ which solve (IV) we must check that $\gamma$
is $d_\mathrm{M}$-exact. Note that
$$d_\mathrm{M} \gamma=\bar{\tau} i_{\derpp{s}}d_s d_\mathrm{M}\Phi_0=0.$$
To prove the exactness it is enough to see that
$$\int_C \gamma=\bar{\tau} i_{\derpp{s}}d_s \int_C \Phi_0=0$$
for every cycle $C$ on $\mathrm{M}$. This holds because $\Phi_0$
represents the same cohomology class on every fibre so $\int_C
\Phi_0$ does not depend on $s$. Therefore there exist 1-forms
$\{\beta_i\}$ on $\mathrm{M}$ depending smoothly on $s$ such that
$\gamma=d_\mathrm{M} \beta_j$. We denote by $\varphi_t$ the
1-parameter group associated to $T$ and by $\mathrm{L}$ the closure
of the abelian group generated by $\varphi_t$. Define
$\tilde{\beta}_j$ as the transverse part of $\int_{\mathrm{L}}
\beta_j(\sigma \cdot)$ (integrating with respect to the Haar measure
on $\mathrm{L}$). Then $\tilde{\beta}_j$ is a transverse 1-form on
$\mathrm{X}$ which depends smoothly on $s$ and such that
$$d_\mathrm{M}\tilde{\beta}_j=d_\mathrm{M}\int_{\mathrm{L}} \sigma^{\ast}\beta_j=
\int_{\mathrm{L}} \sigma^{\ast}
d_\mathrm{M}\beta_j=\int_{\mathrm{L}} \sigma^{\ast}\gamma=\gamma$$
(recall that $\gamma$ is basic, so $\int_{\mathrm{L}}
\sigma^{\ast}\gamma=\gamma$). Therefore there exist basic local
solutions $\{H_i\}$ to (IV) which depends smoothly on $s$. Using a
unit partition as above we obtain a global solution $H_0=\sum_i
\rho(s)H_i$. Note that by construction of $H$ it verifies $\tau
d_\mathrm{M} H=\bar{\tau} d_\mathrm{M} \overline{H}$, therefore it
is also a solution of (II). Now, to proceed with our plan we should
define $F$ as the solution of (III) for the previous $H_0$. Instead
of finding a solution of (III) we will solve the following equation:
$$ d_\mathrm{M} F= \mathrm{i} \big(\tau i_{\derpp{s}}d_s H-\bar{\tau} i_{\derpp{s}}d_s\overline{H} \big)=\nu
\qquad \mathrm{(V)}.$$ Note that the term on the right $\nu$ is a
basic real form on $\mathrm{M}$ depending smoothly on $s$ which is
the sum of a form of type $(1,0)$ and a form of type $(0,1)$,
therefore we can try to apply lemma \ref{lem:Blanch4}. Moreover
$\nu$ is $d_\mathrm{M}$-closed, for
$$d_\mathrm{M} \nu= \mathrm{i} \cdot i_{\derpp{s}}d_s  (\tau d_\mathrm{M} H- \bar{\tau} d_\mathrm{M} \overline{H})=0$$
Nevertheless, here we encounter a difficulty, since we cannot prove
that $\nu_0$, for the previous solution $H_0$, is
$d_\mathrm{M}$-exact. To overcome it we will modify $H_0$ to obtain
a new solution of (IV) for which the corresponding $\nu$ in (V) is
$d_\mathrm{M}$-exact. Consider the basis
$\{\alpha_1,...,\alpha_k,\overline{\alpha}_1,...,\overline{\alpha}_k\}$
of $H^{1,0}(\mathrm{M}/{\mathcal F})\oplus
H^{0,1}(\mathrm{M}/{\mathcal F})$ of forms fixed by $f$ so that they
are well-defined closed forms on $\mathrm{X}$ (note that
$\alpha_1,...,\alpha_k$ are holomorphic). Let
$\{\gamma_1,...,\gamma_{k+1},\overline{\gamma}_1,...,\overline{\gamma}_{k+1}\}$
be the basis of $H_1(\mathrm{M},\C)$ dual to
$\{\alpha_1,...,\alpha_k,\alpha,
\overline{\alpha}_1,...,\overline{\alpha}_k,\overline{\alpha}\}$. We
define
$$u_j=\mathrm{i} \int_{\gamma_j} \big(\tau d_s H-\bar{\tau} d_s\overline{H} \big)
=\mathrm{i} \cdot d_s \int_{\gamma_j} (\tau H-\bar{\tau}
\overline{H})$$ for $1\leq j \leq k$. Notice that
$$\bar{u}_j=\mathrm{i} \cdot d_s \int_{\overline{\gamma}_j} (\tau H-\bar{\tau}
\overline{H})$$ It is not difficult to see that they are
$d_s$-closed and exact real 1-forms on $S^1$ (the base of the
natural fibration). There exists a family $\{v_1,...,v_k\}$ of real
functions on $S^1$ (that we extend by pullback to $\mathrm{X}$) such
that $d_s v_j=u_j$, namely
$$v_j=\mathrm{i} \int_{\gamma_j} (\tau H-\bar{\tau} \overline{H}).$$
We can define now a new solution $H$ of (IV), and therefore a
solution of (I) and (II), by the formula
$$H=H_0+\mathrm{i} \cdot\tau^{-1} \sum_{j=1}^k \alpha_j \cdot v_j$$
so that
$$\nu ds=\mathrm{i}\big(\tau d_s H_0-\bar{\tau} d_s \overline{H}_0\big)+\sum_{j=1}^k u_j\wedge \alpha_j-\sum_{j=1}^k \bar{u}_j\wedge \overline{\alpha}_j.$$
We will next verify that the integral of $\nu ds=\mathrm{i}
\big(\tau d_s H-\bar{\tau} d_s\overline{H})$ is zero for any cycle
$C$ on $\mathrm{M}$, for it is equivalent to $\int_C\nu=0$. Note
that since we saw that $\nu$ is $d_\mathrm{M}$-closed when $H$ is a
solution of (IV) it is enough to check that $\int_{\gamma_j} \nu
ds=\int_{\bar{\gamma}_j} \nu ds=0$ (note that if $\nu=d_\mathrm{M}G$
the function $G$ must be basic). Indeed,
\begin{eqnarray*}
\int_{\gamma_j} \nu ds&=&\int_{\gamma_j} \nu_0 ds-d_s v_j=u_j-u_j=0 \\
\int_{\bar{\gamma}_j} \nu ds&=&\int_{\bar{\gamma}_j} \nu_0
ds-d_s\bar{v}_j=\bar{u}_j-\bar{u}_j=0.
\end{eqnarray*}
Therefore, we are left to solve $d_\mathrm{M} F=\nu$ with $\nu$
satisfying all the hypothesis in lemma \ref{lem:Blanch3} to obtain
local real functions $F_i$ such that $d_\mathrm{M} F_i =\nu$. We
define thus $F$ by means of a unit partition, $F=\sum_i \rho_i F_i$.
To finish the proof it is enough to add a positive constant to $F$
so that we obtain a positive closed $(1,1)$-form on $\mathrm{X}$.
\end{proof}

\subsection{Double suspensions}

We will finally see that in the case of a double suspension we have
a more explicit version of theorem \ref{thm:kahlersusp2} to decide
when the resulting complex manifold is K\"{a}hlerian. If
$\mathrm{N}$ is a compact complex manifold we denote by
Aut$_0(\mathrm{N})$ the connected component of the identity in the
group of holomorphic transformations Aut$_{\C}(\mathrm{N})$ of
$\mathrm{N}$.

\begin{thm}\label{teo:kahlerC2}
Let $\mathrm{N}$ be a compact complex manifold and $f,g
\in\mathrm{Aut}_{\C}(\mathrm{N})$ such that $f\circ g= g\circ f$.
Let $\mathrm{X}$ be the suspension $\mathrm{N} \times \C / \langle
F,G \rangle$ where $F(x,z)=(f(x),z+1)$, $G(x,z)=(g(x),z+\tau)$ and
$\tau\in \C\backslash \R$. Then the following conditions are
equivalent:
\begin{enumerate}[\bf (i)]
\item $\mathrm{X}$ is K\"{a}hlerian.
\item There is a K\"{a}hler form $\omega$ on $\mathrm{N}$ such that
$[f^{\ast}\omega]=[g^{\ast}\omega]=[\omega]$.
\item $\mathrm{N}$ is K\"{a}hlerian and there are integers $n,m>0$ such that
$f^n,g^m\in\mathrm{Aut}_0(\mathrm{N})$.
\end{enumerate}
\end{thm}

\begin{proof}
${\bf (i)\Longrightarrow (ii)}$: Let $\Psi$ be a K\"{a}hler form on
$\mathrm{X}$. Its pull-back $\Phi=\pi^{\ast}\Psi$ by the covering
map $\pi:\mathrm{N}\times \C \rightarrow \mathrm{X}$ is a K\"{a}hler
form on $\mathrm{N}\times \C$. Let us choose $z_0\in \C$. If we
denote by $\mathrm{N}_z=\mathrm{N}\times\{z\}$ and by
$\omega_{z}(x)=\Phi(x,z)_{|\mathrm{N}_z}$ then $\omega_{z_0}$ is a
K\"{a}hler form on $\mathrm{N}_{z_0}$. It suffices to show that
$[f^{\ast}\omega_{z_0}]=[\omega_{z_0}]$ (for $g$ the argument is
analogous). By construction of $\Phi$ we know that
$F^{\ast}\Phi=\Phi$ where $F(x,z)=(f(x),z+1)$, therefore
$f^{\ast}\omega_{z_0}=\omega_{z_0-1}$. Recall that if
$[\Phi]=[\Phi']\in H^2(\mathrm{N}\times \C,\R)$ then
$[\Phi(x,z_0)_{|\mathrm{N}}]=[\Phi'(x,z_0)_{|\mathrm{N}}]$ in
$H^2(\mathrm{N},\R)$. Therefore it is enough to see that
$[\Phi(x,z)]=[\Phi(x,z+a_0)]$ in $H^2(\mathrm{N}\times \C,\R)$ for
all $a_0\in \C$. Since the map $(x,z)\mapsto (x,z+a_0)$ is homotopic
to the identity on $\mathrm{N}\times \C$ this condition is verified
and we can conclude.

${\bf (ii)\Longrightarrow (iii)}$: It is an immediate consequence of
a result by D.Lieberman (c.f. \cite{Lieb}) that asserts the
following: if $\mathrm{M}$ is a compact K\"{a}hler manifold,
$\omega$ a K\"{a}hler form and we denote by
Aut$_{\omega}(\mathrm{M})$ the group of automorphisms of
$\mathrm{M}$ preserving the K\"{a}hler class $[\omega]$ then
Aut$_{\omega}(\mathrm{M})/$Aut$_0(\mathrm{M})$ is a finite group.
Indeed, if we consider $\{f,f^2,...\}$ there must exist $n_1> n_2>0$
such that $f^{n_1}=f^{n_2}\cdot h$ with $h\in$ Aut$_0(\mathrm{N})$.
Therefore for $n=n_1-n_2$ we have $f^n\in$ Aut$_0(\mathrm{N})$ (and
we would proceed identically for $g$).

${\bf (iii)\Longrightarrow (i)}$: It is enough to show that a double suspension of a
compact K\"{a}hler manifold $\mathrm{N}$ by $f,g\in$Aut$_0(\mathrm{N})$ is
K\"{a}hlerian. Indeed, a finite covering $\widetilde{\mathrm{X}}$ of a compact
complex manifold $\mathrm{X}$ is K\"{a}hlerian if and only if $\mathrm{X}$ is
K\"{a}hlerian and the suspension $\mathrm{N}\times \C/ \langle
\widetilde{F},\widetilde{G}\rangle$, for $\widetilde{F}(x,z)=(f^n(x),z+1)$ and
$\widetilde{G}(x,z)=(g^m(x),z+\frac{m}{n}\tau)$, is a finite covering of
$\mathrm{X}$. The statement follows now from a theorem by A.Blanchard (cf.
\cite{Blanch}, p.192) which states the total space $\mathrm{X}$ of a fibred space
with base $\mathrm{B}$ and fibre $\mathrm{F}$ such that $\pi_1(\mathrm{B})$ acts
trivially on $H^1(\mathrm{F},\R)$ is K\"ahlerian if and only if there is a
K\"{a}hler form on $\mathrm{F}$ which represents a cohomology class invariant by
$\pi_1(\mathrm{B})$, $\mathrm{B}$ is K\"{a}hlerian and
$b_1(\mathrm{X})=b_1(\mathrm{B})+b_1(\mathrm{F})$. Note that since
$f,g\in$Aut$_0(\mathrm{N})$ the fundamental group $\pi_1(\mathrm{E})$ acts trivially
on $H^1(\mathrm{N},\R)$ and consequently $b_1(\mathrm{X})=b_1(\mathrm{N})+2$.
\end{proof}

In particular if $\mathrm{N}$ is a compact K\"{a}hler manifold and
Aut$_{\C}(\mathrm{N})$ is finite then any complex manifold
$\mathrm{X}$ thus obtained from $\mathrm{N}$ must be K\"{a}hlerian.

\begin{cor}\label{cor:b1}
Let $\mathrm{N}$ be a compact K\"{a}hler manifold, if
$b_2(\mathrm{N})=1$ any compact complex manifold $\mathrm{X}$
obtained by a double suspension from $\mathrm{N}$ is K\"{a}hlerian.
\end{cor}

\begin{proof}
Applying Hodge decomposition theorem to $\mathrm{N}$ one concludes
that $H^2(\mathrm{N},\C)\cong H^{1,1}(\mathrm{N},\C)\cong \C$. Then
if $\omega$ is a K\"{a}hler form on $\mathrm{N}$ and $f\in
\mathrm{Aut}_{\C}(\mathrm{N})$ we have $[f^{\ast}\omega]=[\lambda
\omega]$, where $\lambda \in \R^+$. Then $f^{\ast}\omega=\lambda
\omega+d\alpha$ where $\alpha$ is a 1-form. Moreover we know that
$[\lambda^n\omega^n]=f^{\ast}[\omega^n]=[\omega^n]$, therefore
$\lambda=\pm 1$. Assume that $\lambda=-1$, then the K\"{a}hler form
$\omega+f^{\ast}\omega=d\alpha$ is exact, which is a contradiction.
Thus $\lambda=1$ and we conclude that the hypothesis $\textbf{(ii)}$
in the theorem \ref{teo:kahlerC2} are trivially fulfilled.
\end{proof}

Although it escapes the scope of this paper note that the same
arguments as in the above theorem yield the following theorem. We
start by recalling the general definition of suspension. Let
$\mathrm{N}$ and $\mathrm{B}$ be two complex manifolds, let
$\widetilde{\mathrm{B}}$ be the universal cover of $\mathrm{B}$ and
$\Gamma=\pi_1(\mathrm{B})$ (acting by the right on
$\widetilde{\mathrm{B}}$). If $h:\Gamma\rightarrow
\mathrm{Aut}_{\C}(\mathrm{N})$ is a group homomorphism we define the
\emph{suspension of $\mathrm{N}$ by $\mathrm{B}$ and $h$} as the
quotient $\Gamma\backslash \mathrm{N}\times \widetilde{\mathrm{B}}$
where the action by $\Gamma$ is
$$\Gamma \times (\mathrm{N}\times\widetilde{\mathrm{B}})\rightarrow (\mathrm{N}\times\widetilde{\mathrm{B}})$$
$$(\gamma,(x,z))\mapsto (h(\gamma)(x), z\cdot \gamma^{-1}).$$

\begin{thm}
Let $\mathrm{N}$ and $\mathrm{B}$ be two compact connected complex
manifolds and $h:\Gamma=\pi_1(\mathrm{B})\rightarrow
\mathrm{Aut}_{\C}(\mathrm{N})$ a group homomorphism. Then the
following conditions are equivalent:
\begin{enumerate}[\bf (i)]
\item The suspension of $\mathrm{N}$ by $\mathrm{B}$ and $h$ is K\"{a}herian.
\item $\mathrm{B}$ is K\"{a}hlerian and there is a K\"{a}hler form $\omega$ on $\mathrm{N}$ such that
$[f^{\ast}\omega]=[\omega]$ for every $f\in h(\Gamma)$. \item
$\mathrm{N}$ and $\mathrm{B}$ are K\"{a}hlerian and there exists a
subgroup $\Gamma'$ of $\Gamma$ such that $h(\Gamma')\subset
\mathrm{Aut}_0(\mathrm{N})$ and $h(\Gamma)/h(\Gamma')$ is finite.
\end{enumerate}
\end{thm}

\bibliographystyle{amsplain}

\providecommand{\bysame}{\leavevmode\hbox
to3em{\hrulefill}\thinspace}
\providecommand{\MR}{\relax\ifhmode\unskip\space\fi MR }
\providecommand{\MRhref}[2]{%
  \href{http://www.ams.org/mathscinet-getitem?mr=#1}{#2}
} \providecommand{\href}[2]{#2}

\end{document}